\documentclass[a4paper,10pt]{article}
\usepackage{amssymb}
\usepackage{latexsym}
\usepackage{amsmath}
\usepackage{amsthm}
\usepackage{amsfonts}
\usepackage{amssymb}
\usepackage{graphicx}

\newtheorem{Th}{Theorem}
\newtheorem{Prop}{Proposition}

\newtheorem{Lm}{Lemma}

\newtheorem{Dfi}{Definition}

\newcommand{\be}{\begin{equation}}
\newcommand{\ee}{\end{equation}}
\newcommand{\R}{\mathbb{R}}

\newcommand{\C}{\mathbb{C}}

\usepackage{mathtools}
\newcommand{\pvint}{\mathop{\mathrlap{\pushpv}}\!\int}
\newcommand{\pushpv}{\mathchoice
  {\mkern5mu\rule[.6ex]{.5em}{1pt}}
  {\mkern2.8mu\rule[.5ex]{.35em}{.8pt}}
  {\mkern2.5mu\rule[.29ex]{.3em}{.7pt}}
  {\mkern2mu\rule[.2ex]{.2em}{.5pt}}
}
\newcommand\res{\mathop{\hbox{\vrule height 7pt width .5pt depth 0pt
\vrule height .5pt width 6pt depth 0pt}}\nolimits}

\def\theequation{\thesection.\arabic{equation}}
\def\theTh{\Roman{section}.\arabic{Th}}
\def\theProp{\Roman{section}.\arabic{Prop}}
\def\theCo{\Roman{section}.\arabic{Co}}

\def\theRm{\Roman{section}.\arabic{Rm}}
\newcommand{\reset}{\setcounter{equation}{0}\setcounter{Th}{0}\setcounter{Prop}{0}\setcounter{Co}{0}
\setcounter{Lm}{0}\setcounter{Rm}{0}}
\def\La{\Lambda}
\def\La{\Lambda}
\def\ti{\tilde}
\def\lf{\left}
\def\rg{\right}

\def\al{\alpha}
\def\la{\lambda}

\def\ep{\varepsilon}
\def\ds{\displaystyle}

\def\Om{\Omega}
\def\om{\omega}
\def\p{\partial}

\def\res{\mathop{\hbox{\vrule height 7pt width .5pt 
depth 0pt\vrule height .5pt width 6pt depth 0pt}}\nolimits}

\begin{document}
\title{The Regularity of Conformal Target Harmonic Maps}
\author{ Tristan Rivi\`ere\footnote{Department of Mathematics, ETH Zentrum,
CH-8093 Z\"urich, Switzerland.}}
\maketitle
{\bf Abstract : }{\it We establish that any weakly conformal $W^{1,2}$ map from
a Riemann surface $S$ into a closed oriented sub-manifold $N^n$ of an euclidian space ${\R}^m$ realizes, for almost every sub-domain, a stationary varifold if and only if it is a smooth 
conformal harmonic map form $S$ into $N^n$.}

\medskip

\noindent{\bf Math. Class. 58E20, 49Q05, 53A10, 49Q15, 49Q20}
\section{Introduction}
In \cite{Ri} the author developed a viscosity method in order to produce closed minimal 2 dimensional surfaces into any arbitrary closed oriented sub-manifolds $N^n$
of any euclidian spaces ${\R}^m$ by min-max type arguments. The method consists in adding to the area of an immersion $\vec{\Phi}$ of a surface $\Sigma$ into $N^n$ a more coercive term such as the $L^{2p}$ norm of the 
second fundamental form preceded by a small parameter $\sigma^2$:
\[
A^\sigma(\vec{\Phi}):=\mbox{Area}(\vec{\Phi})+\sigma^2\int_{\Sigma}\lf[1+|\vec{\mathbb I}_{\vec{\Phi}}|^2\rg]^{p}\ dvol_{g_{\vec{\Phi}}}
\]
where $\vec{\mathbb I}_{\vec{\Phi}}$ is the second fundamental form of the immersion $\vec{\Phi}$ and $dvol_{g_{\vec{\Phi}}}$ is the volume form
associated to the induced metric. For $p>1$ and $\sigma>0$ one proves that the Lagrangians $A^\sigma$ are Palais-Smale in some ad-hoc Finsler bundle of immersions complete for the Palais distance. By applying the now classical {\it Palais-Smale deformation theory} in infinite dimensional space one can then produce critical points $\vec{\Phi}_\sigma$ to $A^\sigma$. It is proved in \cite{MR} that, for a sequence of parameters $\sigma_k\rightarrow 0$, the sequence of integer rectifiable varifolds associated to the immersion of $\Sigma$ by $\vec{\Phi}_{\sigma_k}$ does not necessarily converge\footnote{even modulo extraction of subsequences and in a weak sense
such as the {\it varifold distance} topology.} to a stationary integer rectifiable varifold. However, by applying {\it Struwe's monotonicity trick} one can always select a sequence $\sigma_j\rightarrow 0$ such that the following additional ``entropy estimate'' holds
\[
\sigma_j\ \int_{\Sigma}\lf[1+|\vec{\mathbb I}_{\vec{\Phi}_{\sigma_j}}|^2\rg]^{p}\ dvol_{g_{\vec{\Phi}_{\sigma_j}}}=o\lf(\frac{1}{\sigma_j\,\log\sigma_j^{-1}}\rg)
\]
Assuming this additional estimate, the main achievement of \cite{Ri} is to prove that the immersion of $\Sigma$ by $\vec{\Phi}_{\sigma_j}$ {\it varifold converges} to a  {\it \underbar{stationary} integer rectifiable varifold} given by the image of a smooth Riemann surface $S$ by a weakly conformal $W^{1,2}$ map $\vec{\Phi}$ into $N^n$ equipped by an integer multiplicity. The main result of the paper is to prove that, when this multiplicity is constant, such a map is smooth and satisfies the harmonic map equation. To state our main result we need two definitions.
\begin{Dfi}
\label{almost every}
A property is said to hold for {\bf almost every smooth domain} in $ \Sigma$, if for any smooth domain $\Om$ and any smooth function $f$ such that $f^{-1}(0)=\p\Om$ and $\nabla f\ne 0$ on $\p\Om$ then for almost every $t$ close enough to zero and regular value for $f$  the property holds for the domain contained in $\Omega$ or containing $\Omega$ and bounded by $f^{-1}(\{t\})$.\hfill $\Box$
\end{Dfi}

Precisely we define the notion of {\it target harmonic map} as follows.
\begin{Dfi}
\label{df-I.1}
Let $(\Sigma,h)$ be a smooth closed Riemann surface equipped with a metric compatible with the complex structure. A map $\vec{\Phi}\in W^{1,2}(\Sigma,N^n)$ is {\it target harmonic} if
for almost every smooth domain  $\Om\subset \Sigma$ and any smooth function $F$ supported in the complement of an open neighborhood\footnote{Observe that for almost every domain $\Omega$ the restriction of $\vec{\Phi}$ to $\p\Om$ is H\"older continuous and $\vec{\Phi}(\p\Om)$ is then closed.} of $\vec{\Phi}(\p \Om)$ we have
\be
\label{reg-II.2}
\int_{\Om} \lf<d(F(\vec{\Phi})), d\vec{\Phi}\rg>_h- F(\vec{\Phi})\ A(\vec{\Phi})(d\vec{\Phi},d\vec{\Phi})_h\ dvol_h=0
\ee
where $h$ is any\footnote{One observe that the condition (\ref{reg-II.2}) is conformally invariant and hence independent of the choice of the metric $h$ within the conformal class given by $\Sigma$.} metric compatible with the chosen conformal structure on $S$ and where $A(\vec{q})(\vec{X},\vec{Y})$ denotes the second fundamental form of $N^n$ at the point $\vec{q}$ and acting on the pair of vectors $(\vec{X},\vec{Y})$ and by an abuse of notation we write
\[
\ A(\vec{\Phi})(d\vec{\Phi},d\vec{\Phi})_h:=\sum_{i,j=1}^2h_{ij}\ A(\vec{\Phi})(\p_{x_i}\vec{\Phi},\p_{x_j}\vec{\Phi})\quad.
\]
\hfill $\Box$
\end{Dfi}
Observe that the main difference with the general definition of being {\it harmonic} is that, for {\it target harmonic} one restricts (\ref{reg-II.2}) to test functions $F$ supported in the target
while for the definition of {\it harmonic}, one requires (\ref{reg-II.2}) to hold \underbar{for any} $W^{1,2}$ test function defined on the domain. Therefore, being {\it harmonic} implies to be {\it target harmonic} and the proof of the reverse is the goal of the present work. For a weakly conformal $W^{1,2}$ map into $N^n$ the condition for $\vec{\Phi}$ to be {\it target harmonic} is equivalent to saying that the mapping of $\Sigma$ in $N^n$ defines a  {\it \underbar{stationary} integer rectifiable varifold} (see proposition~\ref{prop-A-1}). Our main result in the present paper is the following.
\begin{Th}
\label{th-I.1}
Any weakly conformal target harmonic map into an arbitrary closed sub-manifold $N^n$ of ${\R}^m$ in  two dimension is smooth and satisfy the harmonic
map equation.\hfill $\Box$
\end{Th}
It is not known if, without the conformality assumption, a {\it target harmonic map} is an harmonic map in the classical sense. Following the main lines of
the proof below one can prove that this is indeed the case in one dimension.

\medskip

{\bf Acknowledgments :}{\it The author would like to thank Alexis Michelat for stimulating discussions on the notion of target harmonic maps and the regularity question for these maps.}
\section{The partial regularity.}

Since the problem is local we shall work in a chart and hence consider maps in $W^{1,2}(D^2,N^n)$ exclusively. Indeed, being {\it target harmonic} in $S$ implies to be {\it target harmonic}
in any sub-domain of $S$. Hence if one proves that {\it target harmonic} on a disc implies {\it harmonic} on that disc in the classical sense we deduce that {\it target harmonic} on $S$ implies {\it harmonic} on $S$. Therefore we can reduce to the case $S=D^2$.

\medskip

Let $\vec{\Phi}$ be a $W^{1,2}$ map from the disc $D^2$ into $N^n$ satisfying the {\it weak conformality condition}
\be
\label{reg-II.1}
|\p_{x_1}\vec{\Phi}|^2=|\p_{x_2}\vec{\Phi}|^2 \quad\mbox{ and }\quad \p_{x_1}\vec{\Phi}\cdot\p_{x_2}\vec{\Phi}=0\quad\mbox{ a.e. on }D^2
\ee
For any $A>1$ we introduce
\[
{\mathcal G}^A:=\lf\{
\begin{array}{c}x\ \mbox{ is a Lebesgue point for }\vec{\Phi}\mbox{ and }\nabla\vec{\Phi}\ \\[3mm] 
x \mbox{ is a point of }L^2\mbox{ app. differentiability }\\[3mm]
 A^{-1}<|\nabla\vec{\Phi}|(x)<A
 \end{array}
 \rg\}
\]
We are going to prove the following lemma
\begin{Lm}
\label{lm-reg-II.1}
Let $\vec{\Phi}$ be a {\it target harmonic map} on $D^2$. Under the previous notations we have that
\[
{\mathcal G}:=\cup_{A>1}{\mathcal G}^A
\]
is an open subset of $D^2$, $\vec{\Phi}$ is smooth and harmonic into $N^n$ in the strong sense on ${\mathcal G}$. Moreover
\be
\label{reg-II.2-a}
{\mathcal H}^2\lf(\vec{\Phi}(D^2\setminus {\mathcal G})\rg)=\frac{1}{2}\int_{D^2\setminus {\mathcal G}}|\nabla\vec{\Phi}|^2(y)\, dy^2=0\quad.
\ee
\hfill $\Box$
\end{Lm}
\noindent{\bf Proof of lemma~\ref{lm-reg-II.1}.}
Let $A>0$. For such an $x\in {\mathcal G}^A$ we shall denote
\[
e^{\la(x)}:=|\p_{x_1}\vec{\Phi}|(x)=|\p_{x_2}\vec{\Phi}|(x)\quad\mbox{ and }\quad \vec{e}_i(x):=e^{-\la(x)}\,\p_{x_i}\vec{\Phi}(x)\quad.
\]
We have in particular 
\[
\frac{1}{|B_r(x)|}\int_{B_r(x)}|\nabla\vec{\Phi}(y)|^2\ dy^2>A^2+o(1)
\]
and hence
\be
\label{reg-II.3}
\lim_{r\rightarrow 0}\frac{\ds\int_{B_r(x)}|\nabla\vec{\Phi}(y)-\nabla\vec{\Phi}(x)|^2\ dy^2}{\ds\int_{B_r(x)}|\nabla\vec{\Phi}(y)|^2\ dy^2}=0\quad.
\ee
And
\be
\label{reg-II.4}
\lim_{r\rightarrow 0}\frac{\ds\int_{B_r(x)}r^{-2}\, |\vec{\Phi}(y)-\vec{\Phi}(x)-\p_{x_1}\vec{\Phi}(x) (y_1-x_1)-\p_{x_2}\vec{\Phi}(x) (y_2-x_2)|^2\ dy^2}{\ds\int_{B_r(x)}|\nabla\vec{\Phi}(y)|^2\ dy^2}=0
\ee
For any $\ep>0$ there exists $r_0>0$ such that for any $r<r_0$ 
\[
\frac{\ds\int_{B_r(x)}|\nabla\vec{\Phi}(y)-\nabla\vec{\Phi}(x)|^2+r^{-2}\, |\vec{\Phi}(y)-\vec{\Phi}(x)-\p_{x_1}\vec{\Phi}(x) (y_1-x_1)-\p_{x_2}\vec{\Phi}(x) (y_2-x_2)|^2\ dy^2}{\ds\int_{B_r(x)}|\nabla\vec{\Phi}(y)|^2\ dy^2}<\ep^2
\]
Using {\it Fubini theorem} together with the {\it mean value Theorem}, for any such $r$ there exists $\rho_{r,x}\in [r/2,r]$ such that
\be
\label{reg-II.4-a}
\begin{array}{l}
\ds\frac{r}{2}\int_{\p B_{\rho_{r,x}(x)}}|\nabla\vec{\Phi}(y)-\nabla\vec{\Phi}(x)|^2+r^{-2}\, |\vec{\Phi}(y)-\vec{\Phi}(x)-\p_{x_1}\vec{\Phi}(x) (y_1-x_1)-\p_{x_2}\vec{\Phi}(x) (y_2-x_2)|^2\ dl_{\p B_{\rho_{r,x}}}\\[5mm]
\ds\quad\quad\quad<\ep^2\int_{B_r(x)}|\nabla\vec{\Phi}(y)|^2\ dy^2= \pi\ r^2\ \ep^2 e^{2\la(x)}(1+o_r(1))
\end{array}
\ee
then
\be
\label{reg-II.5}
\begin{array}{l}
\|\vec{\Phi}(\rho_{r,x},\theta)-\vec{\Phi}(\rho_{r,x},0)-e^{\la(x)}\ \rho_{r,x}\ [\cos\theta-1]\ \vec{e}_1- e^{\la(x)}\ \rho_{r,x}\  \sin\theta\ \vec{e}_2  \|_{L^\infty([0,2\pi])}\\[5mm]
\ds\quad\quad\le \int_0^{2\pi}\lf|\p_\theta\lf( \vec{\Phi}(\rho_{r,x},\theta)-\vec{\Phi}(\rho_{r,x},0)-e^{\la(x)}\  \rho_{r,x}\  [\cos\theta-1]\ \vec{e}_1- e^{\la(x)}\ \rho_{r,x}\ \sin\theta\ \vec{e}_2  \rg)\rg|\ d\theta \\[5mm]
\ds\quad\quad =\int_{\p B_{\rho_{r,x}}(x)} \lf|\frac{1}{\rho_{r,x}}\frac{\p\vec{\Phi}}{\p\theta}(y)+\p_{x_1}\vec{\Phi}(x)\ \sin\theta-\p_{x_2}\vec{\Phi}(x)\ \cos\theta\rg|\ dl_{\p B_{\rho_{r,x}}}\\[5mm]
\ds\quad\quad\le\int_{\p B_{\rho_{r,x}}(x)}|\nabla\vec{\Phi}(y)-\nabla\vec{\Phi}(x)|\ dl_{\p B_{\rho_{r,x}}}<2\,\pi\,\ep\ \rho_{r,x}\ e^{\la(x)} (1+o_r(1))
\end{array}
\ee
Denote $\vec{L}_x(y):=\vec{\Phi}(x)+\p_{x_1}\vec{\Phi}(x) (y_1-x_1)+\p_{x_2}\vec{\Phi}(x) (y_2-x_2)$. Considering (\ref{reg-II.4-a}), we have also chosen $\rho_{r,x}$
in such a way that 
\be
\label{reg-II.6}
\begin{array}{l}
\ds\lf|\pvint_{\p B_{\rho_{r,x}}(x)}\vec{\Phi}(y)\ dl_{\p B_{\rho_{r,x}}}-\vec{\Phi}(x)\rg|^2=\lf|\pvint_{\p B_{\rho_{r,x}}(x)}\vec{\Phi}(y)\ dl_{\p B_{\rho_{r,x}}}-\pvint_{\p B_{\rho_{r,x}}(x)}\vec{L}_x(y)\ dl_{\p B_{\rho_{r,x}}}\rg|^2\\[5mm]
\ds\quad\quad\le \pvint_{\p B_{\rho_{r,x}}(x)}|\vec{\Phi}(y)-\vec{L}_x(y)|^2\ dl_{\p B_{\rho_{r,x}}}\le 2\, \ \ep^2\ \rho^2_{r,x}\ e^{2\la(x)} (1+o_r(1))
\end{array}
\ee
We have 
\be
\label{reg-II.6-a}
\pvint_{\p B_{\rho_{r,x}}(x)}\vec{\Phi}(y)\ dl_{\p B_{\rho_{r,x}}}=\frac{1}{2\pi}\int_0^{2\pi}\vec{\Phi}(\rho_{r,x},\theta)\ d\theta
\ee
We have moreover
\be
\label{reg-II.6-b}
\frac{1}{2\pi}\int_0^{2\pi}\lf[\vec{\Phi}(\rho_{r,x},0)+e^{\la(x)}\ \rho_{r,x}\ [\cos\theta-1]\ \vec{e}_1+e^{\la(x)}\ \rho_{r,x}\  \sin\theta\ \vec{e}_2\rg]\ d\theta=\vec{\Phi}(\rho_{r,x},0)-e^{\la(x)}\ \rho_{r,x}\ \ \vec{e}_1
\ee
Thus combining (\ref{reg-II.5}), ...,(\ref{reg-II.6-b}) we have
\be
\label{reg-II.6.c}
|\vec{\Phi}(x)-\vec{\Phi}(\rho_{r,x},0)+e^{\la(x)}\ \rho_{r,x}\ \ \vec{e}_1|\le \ [2\pi+\sqrt{2}]\ \ep\ \rho_{r,x}\ e^{\la(x)}\ (1+o_r(1))\quad.
\ee
Hence combining (\ref{reg-II.5}) and (\ref{reg-II.6.c}) we obtain
\be
\label{reg-II.7}
\|\vec{\Phi}(\rho_{r,x},\theta)-\vec{\Phi}(x)-e^{\la(x)}\ \rho_{r,x}\ \cos\theta\ \vec{e}_1- e^{\la(x)}\ \rho_{r,x}\  \sin\theta\ \vec{e}_2  \|_{L^\infty([0,2\pi])}\le 5\,\pi\, \ep\ \rho_{r,x}\ e^{\la(x)}\ (1+o_r(1))
\ee
Denote by $\Gamma_{\rho_{r,x}}$ the circle in ${\R}^m$ given by $\vec{L}_x(\p B_{\rho_{r,x}}(x))$ and
\be
\label{reg-II.8}
\Gamma_{\rho_{r,x}}^\ep=\lf\{y\in {N^n}\quad\mbox{ s. t. }\quad\mbox{dist}(y,\Gamma_{\rho_{r,x}})>6\,\pi\, \ep\ \rho_{r,x}\ e^{\la(x)}\rg\}
\ee
Because of our assumptions,  $\Sigma_{\rho_{r,x}}:=\vec{\Phi}(B_{\rho_{r,x}}(x))\cap \Gamma_{\rho_{r,x}}^\ep$ defines a rectifiable integer stationary varifold in $N^n\cap \Gamma_{\rho_{r,x}}^\ep$. Leon Simon monotonicity formula in ${\R}^m$ gives 
\be
\label{reg-II.7-a}
\begin{array}{l}
\ds\rho_{r,x}^{-2}\ {\mathcal H}^2\lf(\Sigma_{\rho_{r,x}}\cap B^m_{e^{\la(x)}\,\rho_{r,x}}(\vec{\Phi}(x))\rg)-\pi\, e^{2\,\la(x)}\,\ge -\frac{1}{4}\int_{B^m_{e^{\la(x)}\,\rho_{r,x}}(\vec{\Phi}(x))}  |\vec{H}_{{\R}^m}|^2 d{\mathcal H}^2\res\Sigma_{\rho_{r,x}}\\[5mm]
\ds\quad\quad-\frac{1}{\rho_{r,x}^{2}}\int_{\vec{q}\in B^m_{e^{\la(x)}\,\rho_{r,x}}(\vec{\Phi}(x))}  |\vec{H}_{{\R}^m}| |\vec{q}-\vec{\Phi}(x)| d{\mathcal H}^2\res\Sigma_{\rho_{r,x}}\end{array}
\ee
where $\vec{H}_{{\R}^m}$ is the generalized curvature of the varifold given by $\vec{\Phi}$ in ${\R}^m$. The generalized curvature is by definition given by
\[
\int \vec{H}_{{\R}^m}\cdot F(\vec{q})\ d{\mathcal H}^2\ =-\,\int_\Sigma d(F(\vec{\Phi}))\cdot d\vec{\Phi}\ dvol_{g_{\vec{\Phi}}}
\]
Using the stationarity in $N^n$ we have (\ref{A-2}) and then we have
\[
\vec{H}_{{\R}^m}:=-\,A(\vec{\Phi})(d\vec{\Phi},d\vec{\Phi})_{g_{\vec{\Phi}}}\quad.
\]
Hence, since $|d\vec{\Phi}|_{g_{\vec{\Phi}}}=1$ we deduce that $|\vec{H}_{{\R}^m}|\le \|A\|_{L^\infty(N^n)}$. Inserting this bound in (\ref{reg-II.7-a}) we obtain
\be
\label{reg-II.7-b}
\rho_{r,x}^{-2}\ {\mathcal H}^2\lf(\Sigma_{\rho_{r,x}}\cap B^m_{e^{\la(x)}\,\rho_{r,x}}(\vec{\Phi}(x))\rg)-\pi\, e^{2\,\la(x)}\,\ge -\,C\, \rho_{r,x}^{-1}\ {\mathcal H}^2\lf(\Sigma_{\rho_{r,x}}\cap B^m_{e^{\la(x)}\,\rho_{r,x}}(\vec{\Phi}(x))\rg)
\ee
Thus
\be
\label{reg-II.7-c}
\rho_{r,x}^{-2}\ {\mathcal H}^2\lf(\Sigma_{\rho_{r,x}}\cap B^m_{e^{\la(x)}\,\rho_{r,x}}(\vec{\Phi}(x))\rg)\ge\pi\, e^{2\,\la(x)}\,(1-C\,\rho_{r,x})\ \ee
Thus having chosen $r_0<\ep^2$, we have
\be
\label{reg-II.8-a}
{\mathcal H}^2\lf(\Sigma_{\rho_{r,x}}\cap B^m_{e^{\la(x)}\,\rho_{r,x}}(\vec{\Phi}(x))\rg)\ge \pi\ e^{2\,\la(x)}\ \rho_{r,x}^2\ \lf(1-C\ep^2\rg)\quad.
\ee
In the mean time, due to (\ref{reg-II.3}) we have
\be
\label{reg-II.9}
{\mathcal H}^2\lf(\Sigma_{\rho_{r,x}}\rg)\le\frac{1}{2}\int_{B_{\rho_{r,x}}(x)}|\nabla\vec{\Phi}(y)|^2\ dy^2\le \pi\ e^{2\la(x)}\ \rho^2_{r,x}\ (1+C\,\ep^2)\quad.
\ee
The {\it tilt excess} of $ \Sigma_{\rho_{r,x}}$ in the sense of Allard \cite{All} is given by
\[
E\lf(\Sigma_{\rho_{r,x}},\vec{e}_1\wedge\vec{e}_2,\vec{\Phi}(x),\rho_{r,x}\rg)= \frac{1}{\rho_{r,x}^2}\int_{B_{\rho_{r,x}}(x)\cap {\mathcal G}}|e^{-2\la(y)}\p_{x_1}\vec{\Phi}(y)\wedge\p_{x_2}\vec{\Phi}(y)-\vec{e}_1\wedge\vec{e}_2|^2\ e^{2\la(y)}\ dy^2
\]
Since we have the pointwize bound 
\be
\label{reg-II.9-a}
2\, |e^{\la(y)}-e^{\la(x)}|=||\nabla\vec{\Phi}|(x)-|\nabla\vec{\Phi}|(y)|\le |\nabla\vec{\Phi}(x)-\nabla\vec{\Phi}(y)|
\ee
hence we deduce from (\ref{reg-II.4-a}) that\footnote{We are restricting to the integration on  $B_{\rho_{r,x}}(x)\cap {\mathcal G}$ because one requires $\la$ to be defined as a measurable function.}
\[
\int_{B_{\rho_{r,x}}(x)\cap {\mathcal G}}|e^{\la(y)}-e^{\la(x)}|^2\ dy^2\le \pi\ \ep^2\ \rho_{r,x}^2 e^{2\la(x)}
\]
Hence, denoting
\[
E_\ep:=\lf\{y\in B_{\rho_{r,x}}(x)\cap {\mathcal G}\ ;\ |\la(y)-\la(x)|> \sqrt{\ep}\rg\}\quad,
\]
we have
\be
\label{reg-II.10}
|E_\ep|\le \ep\ |B_{\rho_{r,x}}(x)|\quad.
\ee
Observe
\be
\label{reg-II.11}
\begin{array}{l}
\ds\int_{E_\ep}|\nabla\vec{\Phi}(y)|^2\ dy^2\le 2\, \int_{B_{\rho_{r,x}}(x)}|\nabla\vec{\Phi}(y)-\nabla\vec{\Phi}(x)|^2\ dy^2+2\, \int_{E_\ep}|\nabla\vec{\Phi}(x)|^2\ dy^2\\[5mm]
\ds\quad\quad\le 4\,\pi\, [\ep+\ep^2]\ \rho_{r,x}^2\,e^{2\la(x)}\le 3\, \pi\, \ep\ \rho_{r,x}^2\,e^{2\la(x)}
\end{array}
\ee
Observe that we have the pointwize bound
\[
\begin{array}{l}
|e^{-2\la(y)}\p_{x_1}\vec{\Phi}(y)\wedge\p_{x_2}\vec{\Phi}(y)-\vec{e}_1\wedge\vec{e}_2|^2\ e^{2\la(y)}=\ |\vec{e}_1(y)\wedge\p_{x_2}\vec{\Phi}-\vec{e}_1(x)\wedge\vec{e}_2(x)\,e^{\la(y)}|^2\\[5mm]
\ds\quad\le2\, |\vec{e}_1(x)-\vec{e}_1(y)|^2\ e^{2\la(y)}+2\, |\vec{e}_2(x)-\vec{e}_2|^2 \ e^{2\la(y)}\le 4\, |e^{\la(x)}-e^{\la(y)}|^2+2\, |\nabla\vec{\Phi}(x)-\nabla\vec{\Phi}(y)|^2\\[5mm]
\ds\quad\le 4\, |\nabla\vec{\Phi}(x)-\nabla\vec{\Phi}(y)|^2
\end{array}
\]
where we have also used (\ref{reg-II.9-a}). We have using one more time (\ref{reg-II.4-a}) together with  (\ref{reg-II.11}) and the previous pointwize inequality
\be
\label{reg-II.12}
\begin{array}{l}
\ds E\lf(\Sigma_{\rho_{r,x}},\vec{e}_1\wedge\vec{e}_2,\vec{\Phi}(x),\rho_{r,x}\rg)\le \frac{4}{\rho_{r,x}^2}\int_{E_\ep}|\nabla\vec{\Phi}(y)|^2\ dy^2\\[5mm]
\ds+\frac{1}{\rho_{r,x}^2}\int_{B_{\rho_{r,x}}(x)\cap {\mathcal G}\setminus E_\ep}|e^{-2\la(y)}\p_{x_1}\vec{\Phi}(y)\wedge\p_{x_2}\vec{\Phi}(y)-\vec{e}_1\wedge\vec{e}_2|^2\ e^{2\la(y)}\ dy^2\\[5mm]
\ds\quad\le 12\, \pi\, \ep\ A^2+\ep\, \frac{1}{\rho_{r,x}^2}\int_{B_{\rho_{r,x}}(x)}|\nabla\vec{\Phi}(y)|^2\ dy^2\le 14\, \pi\, \ep\ A^2\quad.
\end{array}
\ee
where we recall that $x$ has been chosen in such a way that $2\,e^{2\la(x)}\le A^2$. Combining (\ref{reg-II.9}) and (\ref{reg-II.12}) we can apply Allard main regularity result and precisely we obtain for $\ep$ small enough the existence of $\gamma<1$ such that ${\mathcal S}_{r,x}:=\Sigma_{\rho_{r,x}}\cap B^m_{\gamma\, e^{\la(x)}\,r}(\vec{\Phi}(x))$ is a smooth minimal sub-manifold. Moreover, because of the {\it upper semi-continuity} of the density function for a stationarity varifold we can assume that all points $\vec{q}\in {\mathcal S}_{r,x}$ have multiplicity 1
(See for instance the presentation of the absence of hole in \cite{Del} section 7). In other words we have that
\be
\label{reg-II.13}
\begin{array}{l}
\ds\int_{B^m_{\gamma\, e^{\la(x)}\,r}(\vec{\Phi}(x))}d{\mathcal H}^2\res\Sigma_{\rho_{r,x}}=\mbox{Area}({\mathcal S}_{r,x}) =\int_{\vec{\Phi}^{-1}({\mathcal S}_{r,x})}dvol_{g_{\vec{\Phi}}}=\frac{1}{2}\int_{\vec{\Phi}^{-1}( {\mathcal S}_{r,x})}|\nabla\vec{\Phi}|^2\ dx^2
\end{array}
\ee
and, in the area formula, the counting function 
\be
\label{counting}
{\mathcal H}^0(\vec{\Phi}^{-1}({\mathcal S}_{r,x})\cap\vec{\Phi}^{-1}\{\vec{q}\})=1\quad\quad\mbox{ for }\quad{\mathcal H}^2-\mbox{ a. e.}\quad\vec{q}\in{\mathcal S}_{r,x}\quad.
\ee

\medskip

We claim that for $s$ small enough $\vec{\Phi}(B_s(x))\subset {\mathcal S}_{r,x}$. Assuming there is point $\vec{p}\not\in B^m_{\gamma r\ e^{\la}(x)}(\vec{\Phi}(x))$ but $\vec{p}\in \vec{\Phi}(B_s(x))$. As $s$ goes to zero we have in one hand, using (\ref{reg-II.4-a}),
\be
\label{reg-II.13}
\int_{B_s(x)}|\nabla\vec{\Phi}|^2\ dy^2\le 2\,\pi\, s^2\, (1+\ep) \, e^{2\la(x)} 
\ee
 but the stationarity of $\Sigma_{\rho_{s,x}}:=\vec{\Phi}(B_{\rho_{s,x}}(x))\cap \Gamma_{\rho_{2\,s,x}}^\ep$ in $N^n\cap \Gamma_{\rho_{2\,s,x}}^\ep$ together with the monotonicity formula would imply
 \[
 \pi (\gamma r-s)^2\, e^{2\la(x)}\ (1-o_r(1))\le \pi (\mbox{dist}(\vec{p},\Gamma_{\rho_{2\,s,x}}^\ep)\ (1-o_r(1)) \le {\mathcal H}^2(\Sigma_{\rho_{2\,s,x}})\le \frac{1}{2}\int_{B_s(x)}|\nabla\vec{\Phi}|^2\ dy^2
 \]
 which contradicts\footnote{Observe that based on similar arguments one could prove directly that $\vec{\Phi}$ is continuous on the whole $\Sigma$ without using Allard's result
 but this is not needed.} (\ref{reg-II.13}) for $s$ very small compared to $r$. Hence, for $s$ small enough $$\vec{\Phi}(B_s(x))\subset {\mathcal S}_{r,x}= \Sigma_{\rho_{r,x}}\cap B^m_{\gamma \, e^{\la(x)} \, r}(\vec{\Phi}(x))$$

 Let $\Psi$ be a smooth conformal diffeomorphism from $\Sigma_{\rho_{r,x}}\cap B^m_{\gamma \, e^{\la(x)} \, r}(\vec{\Phi}(x))$ into
the disc $D^2$. Denote by $e^{2\mu}\ dz^2=\Psi^{\ast} g_{{\R}^m}$ and consider $f:=\Psi\circ\vec{\Phi}$. This is in particular a $W^{1,2}$ weakly conformal map from $B_s(x)$ into $\C$. Consider any map $g\in W^{1,2}(B_{\rho_{s,x}}(x))$ such that $g=f$ on $\p B_{\rho_{s,x}}(x)$. The union of $f$ and $g$ realizes a $W^{1,2}\cap L^\infty$ map from $S^2$ into ${\C}$ that we denote $\ti{f}$. The cycle $\ti{f}_\ast[S^2]$ has then an {\it algebraic covering number} equal to $0$ ${\mathcal H}^2-$almost everywhere in ${\C}$ that is to say\footnote{Indeed due to the $W^{1,2}$ nature of $\ti{f}$ and since we are in 2 dimension we
have that $d$ and $\ast$ commute (this is clearly not the case in higher dimension) and we have in particular $\ti{f}^\ast dz_1\wedge dz_2=d(\ti{f}_1\,d\ti{f}_2)$}
\[
\int_{S^2}\ti{f}^\ast dz_1\wedge dz_2=0\quad.
\] 
Away from the compact 1 rectifiable set $f(\p B_{\rho_{s,x}}(x))$, because of (\ref{counting}), the image $f(B_{\rho_{s,x}}(x))$ has covering number $+1$, $-1$ or $0$. Hence, $g$ must have an odd covering number almost everywhere in ${\C}\setminus f(\p B_{\rho_{s,x}}(x))$ whenever the covering number of $f$ is non zero. This homological fact implies
\[
\begin{array}{l}
\ds\frac{1}{2}\int_{B_{\rho_{s,x}}(x))}e^{2\mu(f)}\,|\nabla f|^2\ dx^2=\int_{B_{\rho_{s,x}}(x))}e^{2\mu(f)}|\p_{x_1}f\times\p_{x_2}f|\ dx_1\wedge dx_2\\[5mm]
\ds\quad\le\int_{B_{\rho_{s,x}}(x))}e^{2\mu(g)}|\p_{x_1}g\times\p_{x_2}g|\ dx_1\wedge dx_2\le\frac{1}{2}\int_{B_{\rho_{s,x}}(x))}e^{2\mu(g)}\,|\nabla g|^2\ dx^2
\end{array}
\]
Hence by the {\it Dirichlet principle} $f$ coincides with it's harmonic extension\footnote{The harmonic extension is unique for r small enough since one is taking value in a convex
geodesic ball of $(D^2,e^{2\mu}\ dz^2)$.}  in $(D^2,e^{2\mu}\ dz^2)$. The map $\Psi\circ\vec{\Phi}$ is then smooth (holomorphic or anti-holomorphic) on $B_{\rho_{s,x}}(x))$. This implies Lemma~\ref{lm-reg-II.1}.\hfill $\Box$


\section{The Harmonicity of Conformal Target Harmonic Maps}

\subsection{The integration by part formula.}

The goal of the present subsection is to prove the following {\it integration by parts formula.}
\begin{Lm}
\label{lm-reg-III.5}
Let $F(\vec{p})=(F_{ij}(\vec{p}))_{{1\le i\le m}\ {1\le j\le k}}$ be a $C^1$ map on $N^n$ taking values into $m\times k$ real matrices. Then, for any 
smooth compactly supported function $\varphi$ in $D^2$ and for almost every regular value $t>0$ of $\varphi$ one has
\be
\label{reg-III.5}
\sum_{i=1}^m\int_{\p \Om_t} F(\vec{\Phi}){ij}\,\frac{\p \vec{\Phi}_i}{\p \nu}\, dl_{\p \Om_t}-\int_{\Om_t} \nabla(F(\vec{\Phi}){ij})\, \nabla\vec{\Phi}_i\ dy^2+\int_{\Om_t}F(\vec{\Phi})_{ij}\, A_i(\vec{\Phi})(d\vec{\Phi},d\vec{\Phi})\ dy^2=0
\ee
where $\Om_t=\varphi^{-1}((t,+\infty))$ and $\nu$ is the exterior unit normal to the level set $\nu:=\nabla\varphi/|\nabla\varphi|$.
\hfill $\Box$
\end{Lm}
In order to prove lemma~\ref{lm-reg-III.5}, we will first establish some intermediate results.

\medskip

Let $x\in D^2$ and choose $r$ such that
\be
\label{reg-III.1}
\frac{1}{2\ep}\int_{r-\ep}^{r+\ep}ds\int_{\p B_s(x)}|\nabla\vec{\Phi}|^2\ dl_{\p B_s}=\int_{\p B_r(x)}|\nabla\vec{\Phi}|^2\ dl_{\p B_r}<+\infty
\ee
Hence the restriction of $\vec{\Phi}$ to $\p B_r(x)$ is $W^{1,2}$ and the continuous image of $\p B_r(x)$ by $\vec{\Phi}$, $\Gamma_{r,x}:=\vec{\Phi}(\p B_r(x))$ has finite length and is rectifiable and compact.

We denote ${\mathcal B}:=D^2\setminus {\mathcal G}$. Because of the previous lemma~\ref{lm-reg-II.1} we have that ${\mathcal B}$ is closed. Hence
\[
{\mathcal B}=\bigcap_{\ep>0} {\mathcal B}_\ep\quad\mbox{where }\quad {\mathcal B}_\ep:=\lf\{x\in D^2\ ;\ \mbox{dist}(x,{\mathcal B})\le\ep\rg\}
\]
Since the integral of $|\nabla\vec{\Phi}|^2$ over ${\mathcal B}$ is zero, we clearly have
\be
\label{reg-III.1-a}
\lim_{\ep\rightarrow 0}\int_{{\mathcal B}_\ep}|\nabla\vec{\Phi}|^2(y)\ dy^2=0
\ee
We denote also
\[
{\mathcal G}_\ep:=D^2\setminus{\mathcal B}_\ep
\]
Finally we denote by ${\frak B}_\ep$ the {\it rectifiable image} of ${\mathcal B}_\ep$ by $\vec{\Phi}$ that is the image by the approximate continuous
representative of $\vec{\Phi}$ of the intersection of ${\mathcal B}_\ep$ with the points of approximate differentiability of $\vec{\Phi}$. We claim the following.
\begin{Lm}
\label{lm-reg-III.2}
Under the previous notations we have
\be
\label{reg-III.3}
\lim_{\ep\rightarrow 0}\int_{D^2\cap \vec{\Phi}^{-1}({\frak B}_\ep)}|\nabla\vec{\Phi}|^2(y)\ dy^2=0
\ee
\end{Lm}
\noindent{\bf Proof of lemma~\ref{lm-reg-III.2}.}  Identity (\ref{reg-III.1-a}) implies that
\be
\label{reg-III.4}
\lim_{\ep\rightarrow 0}{\mathcal H}^2({\frak B}_\ep)=0
\ee
Hence for any $\delta>0$ there exists $\ep_\delta>0$ and for any $\ep<\ep_\delta$ there exists a covering of ${\frak B}_\ep$ by balls $(B^m_{r_l}(\vec{p}_l))_{l\in \Lambda}$ such that
\[
\sum_{l\in \Lambda} r_l^2<\delta
\]
Using the monotonicity formula for stationary varifolds we obtain
\[
\int_{D^2\cap \vec{\Phi}^{-1}({\frak B}_\ep)}|\nabla\vec{\Phi}|^2(y)\ dy^2\le\sum_{i\in I}\int_{\vec{\Phi}^{-1}(B^m_{r_l}(\vec{p}_l))}|\nabla\vec{\Phi}|^2(y)\ dy^2\le C\, \sum_{l\in \Lambda} r_l^2<C\,\delta
\]
which implies the lemma.\hfill $\Box$

\medskip

We shall prove the following.
\begin{Lm}
\label{lm-reg-III.3}
Under the previous notations we have
\be
\label{reg-III.4-1}
\lim_{s\rightarrow 0}\int_{{\mathcal G}_\ep}\lf[1+|\nabla\vec{\Phi}|^2(y)\rg]\lf|\frac{{\mathcal H}^2\lf({\mathfrak B}_\ep\cap B^m_s(\vec{\Phi}(y))\rg)}{\pi\ s^2}-\theta^2_\ep(\vec{\Phi}(y))\rg|\ dy^2=0
\ee
where
\[
\theta^2_\ep(\vec{p}):=\lim_{s\rightarrow 0}\frac{{\mathcal H}^2\lf({\mathfrak B}_\ep\cap B^m_s(\vec{p})\rg)}{\pi\ s^2}
\]
exists ${\mathcal H}^2$ almost everywhere since ${\mathfrak B}_\ep\subset\cup_{\eta>0}\vec{\Phi}({\mathcal G}_\eta)$ is 2-rectifiable. Moreover
since $\vec{\Phi}$ is a smooth immersion on ${\mathcal G}_\ep$, $\theta^2_\ep(\vec{\Phi}(y))$ is a well defined measurable function on ${\mathcal G}_\ep$.
\end{Lm}
\noindent{\bf Proof of lemma~\ref{lm-reg-III.3}.} Since $|\nabla\vec{\Phi}|$ is uniformly bounded on ${\mathcal G}_\ep$ we have
\[
\begin{array}{l}
\ds\int_{{\mathcal G}_\ep}\lf[1+|\nabla\vec{\Phi}|^2(y)\rg]\lf|\frac{{\mathcal H}^2\lf({\mathfrak B}_\ep\cap B^m_s(\vec{\Phi}(y))\rg)}{\pi\ s^2}-\theta^2_\ep(\vec{\Phi}(y))\rg|\ dy^2\\[5mm]
\ds\quad\le \lf[1+\|\nabla\vec{\Phi}\|_{L^\infty({\mathcal G}_\ep)}\rg]\int_{{\mathcal G}_\ep}\lf|\frac{{\mathcal H}^2\lf({\mathfrak B}_\ep\cap B^m_s(\vec{\Phi}(y))\rg)}{\pi\ s^2}-\theta^2_\ep(\vec{\Phi}(y))\rg|\ dy^2
\end{array}
\]
Since  $\vec{\Phi}$ is a smooth immersion on ${\mathcal G}_\ep$, $\theta^2_\ep(\vec{\Phi}(y))$ is a well defined measurable function on ${\mathcal G}_\ep$ and we have for almost $y$
\[
\lim_{s\rightarrow 0}\frac{{\mathcal H}^2\lf({\mathfrak B}_\ep\cap B^m_s(\vec{\Phi}(y))\rg)}{\pi\ s^2}=\theta^2_\ep(\vec{\Phi}(y))
\]
The monotonicity formula for the stationary varifold given by the image of $D^2$ by $\vec{\Phi}$ gives
\[
\sup_{s>0;\ y\in D^2}\frac{{\mathcal H}^2\lf({\mathfrak B}_\ep\cap B^m_s(\vec{\Phi}(y))\rg)}{\pi\ s^2}\le \sup_{s>0;\ y\in D^2}\frac{{\mathcal H}^2\lf(\vec{\Phi}(D^2)\cap B^m_s(\vec{\Phi}(y))\rg)}{\pi\ s^2}\le C
\]
for some $C$. The lemma follows by a direct application of dominated convergence.\hfill $\Box$

\medskip

\noindent{\bf Proof of lemma~\ref{lm-reg-III.5}} In order to simplify the presentation we shall restrict to $\Omega_t$ to be balls. Because of (\ref{reg-III.3}), for every $x\in D^2$ and almost every $r>0$ we have
\be
\label{reg-III.6}
\lim_{\ep\rightarrow 0}\int_{\vec{\Phi}^{-1}({\mathfrak B}_\ep)\cap \p B_r(x)}|\nabla\vec{\Phi}|(y)\ dl_{\p B_r}=0
\ee.
We choose $r$ such that (\ref{reg-III.6}) holds true and such that also
\be
\label{reg-III.5-a}
\lim_{s\rightarrow 0}\int_{{\mathcal G}_\ep\cap\p B_r(x)}\lf[1+|\nabla\vec{\Phi}|^2(y)\rg]\lf|\frac{{\mathcal H}^2\lf({\mathfrak B}_\ep\cap B^m_s(\vec{\Phi}(y))\rg)}{\pi\ s^2}-\theta^2_\ep(\vec{\Phi}(y))\rg|\ dl_{\p B_r}=0\quad.
\ee
Hence in particular ${\mathcal H}^1(\vec{\Phi}({\mathcal B}_\ep\cap \p B_r(x)))$ is converging to zero as $\ep$ goes to zero. Since $\vec{\Phi}$ is continuous on $\p B_r(x)$ and since ${\mathcal B}_\ep$ is a closed set we have that $\vec{\Phi}({\mathcal B}_\ep\cap \p B_r(x))$ is a compact subset
of $N^n$. Because of the previous, for any $\delta>0$, there exists $\ep_\delta>0$ such that for any $\ep<\ep_\delta$, we can include $\vec{\Phi}({\mathcal B}_\ep\cap \p B_r(x))$ in \underbar{finitely}
many balls $(B^m_{r_l}(\vec{p}_l))_{l\in\Lambda}$ such that
\be
\label{reg-III.7}
\sum_{l\in \Lambda}r_l<\delta\quad.
\ee
We also choose $\ep_\delta>0$ such that for any $\ep<\ep_\delta$
\be
\label{reg-III.7-a}
\int_{\vec{\Phi}^{-1}({\mathfrak B}_\ep)\cap \p B_r(x)}|\nabla\vec{\Phi}|(y)\ dl_{\p B_r}<\delta\quad.
\ee
Let $\chi$ be a cut-off function on ${\R}_+$ such that $\chi\equiv 1$ on $[2,+\infty)$ and $\chi\equiv 0$ on $[0,1]$. We introduce
\[
\xi_\ep(\vec{\Phi}(y)):=\prod_{l\in \Lambda}\chi\lf(\frac{|\vec{\Phi}(y)- \vec{p}_l|}{r_l}\rg)
\]
Observe that $\xi_\ep(\vec{\Phi}(y))$ is zero on ${\mathcal B}_\ep\cap \p B_r(x)$. For any $s>0$ we also introduce
\[
\eta_s(\vec{\Phi}(y)):=\chi\lf(\frac{\mbox{dist}(\vec{\Phi}(y),\Gamma_{r,x})}{s}\rg)
\]
where we recall that $\Gamma_{r,x}:=\vec{\Phi}(\p B_r(x))$. Using the assumption (\ref{reg-II.2}) we have for any $F$ as in the statement 
of the lemma
\be
\label{reg-III.8}
\sum_{i=1}^m-\int_{B_r(x)} \nabla(\xi_\ep(\vec{\Phi})\, \eta_s(\vec{\Phi})\,F(\vec{\Phi})_{ij})\, \nabla\vec{\Phi}_i\ dy^2+\int_{B_r(x)}\xi_\ep(\vec{\Phi})\, \eta_s(\vec{\Phi}) F(\vec{\Phi})_{ij}\, A_i(\vec{\Phi})(d\vec{\Phi},d\vec{\Phi})\ dy^2=0
\ee
We have
\be
\label{reg-III.9}
\begin{array}{l}
\ds\lf|\int_{B_r(x)} \nabla(\xi_\ep(\vec{\Phi}))\, \eta_s(\vec{\Phi})\,F(\vec{\Phi})_{ij}\, \nabla\vec{\Phi}_i\ dy^2\rg|\\[5mm]
\ds\le C\ \|F\|_\infty\,\sum_{l\in \La}\frac{1}{r_l}\int_{B_r(x)}|\chi'|\lf(\frac{|\vec{\Phi}(y)- \vec{p}_l|}{r_l}\rg)|\nabla\vec{\Phi}|^2(y)\ dy^2\\[5mm]
\ds\le C\ \|F\|_\infty\, \sum_{l\in \Lambda} r_l^{-1}\ \int_{\vec{\Phi}^{-1}(B^m_{r_l}(\vec{p}_l))}|\nabla\vec{\Phi}|^2(y)\ dy^2\le\
 C\ \|F\|_\infty\, \sum_{l\in \Lambda}r_l\le C\ \|F\|_\infty\, \delta
\end{array}
\ee
where we observe that the bound is independent of $s$. We now write
\be
\label{reg-III.10}
\int_{B_r(x)} \xi_\ep(\vec{\Phi})\, \nabla(\eta_s(\vec{\Phi}))\,F(\vec{\Phi})_{ij}\, \nabla\vec{\Phi}_i\ dy^2=\int_{B_r(x)\cap {\mathcal B}_\ep}\cdots+\int_{B_r(x)\cap {\mathcal G}_\ep}\cdots
\ee
We have
\be
\label{reg-III.11}
\begin{array}{l}
\ds\lf|\int_{B_r(x)\cap{\mathcal B}_\ep} \xi_\ep(\vec{\Phi})\, \nabla(\eta_s(\vec{\Phi}))\,F(\vec{\Phi})_{ij}\, \nabla\vec{\Phi}_i\ dy^2\rg| \\[5mm]
\ds\quad\le \frac{C\ \|F\|_\infty}{s}\int_{B_r(x)\cap{\mathcal B}_\ep} {\mathbf 1}_{\mbox{dist}(\vec{\Phi}(y),\Gamma_{r,x}^\ep)<s}\ |\nabla\vec{\Phi}|^2(y)\ dy^2
\end{array}
\ee
Where $\Gamma_{r,x}^\ep$ is the smooth immersed curve $\vec{\Phi}(\p B_r(x))\cap{\mathcal G}_\ep$ and ${\mathbf 1}_{\mbox{dist}(\vec{\Phi}(y),\Gamma_{r,x}^\ep)<s}(y)$ is the characteristic function of the set of $y$ such that $\vec{\Phi}(y)$ is at the distance at most $s$ to $\Gamma_{r,x}^\ep$.
The fact that we can restrict to $\Gamma_{r,x}^\ep$ instead of $\Gamma_{r,x}$ is due to the fact that we are cutting off $\Gamma_{r,x}\setminus\Gamma_{r,x}^\ep$ by multiplying by $\xi_\ep(\vec{\Phi})$. Observe that since the curve $\Gamma_{r,x}^\ep$ is a smooth immersion of the open subset
of $\p B_r(x)$ given by $\p B_r(x)\cap {\mathcal G}_\ep$ we have for $s$ small enough
\be
\label{reg-III.12}
\begin{array}{l}
\ds{\mathbf 1}_{\mbox{dist}(\vec{p},\Gamma_{r,x}^\ep)<s}\le \frac{1}{s}\int_{\Gamma_{r,x}^\ep}{\mathbf 1}_{\mbox{dist}(\vec{p},\vec{q})<2s}\ d{\mathcal H}^1(\vec{q})\\[5mm]
\ds\quad\le  \frac{1}{s}\int_{\p B_r(x)\cap {\mathcal G}_\ep}{\mathbf 1}_{\mbox{dist}(\vec{p},\vec{\Phi}(z))<2s}\ |\nabla\vec{\Phi}|(z)\ dl_{\p B_r}
\end{array}
\ee
Inserting this inequality in (\ref{reg-III.11}) gives
\be
\label{reg-III.13}
\begin{array}{l}
\ds\lf|\int_{B_r(x)\cap{\mathcal B}_\ep} \xi_\ep(\vec{\Phi})\, \nabla(\eta_s(\vec{\Phi}))\,F(\vec{\Phi})_{ij}\, \nabla\vec{\Phi}_i\ dy^2\rg| \\[5mm]
\ds\quad\le \frac{C\ \|F\|_\infty}{s^2}\int_{\p B_r(x)\cap {\mathcal G}_\ep}\ |\nabla\vec{\Phi}|(z)\ dl_{\p B_r}\int_{B_r(x)\cap{\mathcal B}_\ep}{\mathbf 1}_{\mbox{dist}(\vec{\Phi}(y),\vec{\Phi}(z))<2s}\ |\nabla\vec{\Phi}|^2(y)\ dy^2\\[5mm]
\ds\quad\le C\ \|F\|_\infty\ \int_{\p B_r(x)\cap {\mathcal G}_\ep}\ |\nabla\vec{\Phi}|(z)\ \frac{{\mathcal H}^2({\mathfrak B}_\ep\cap B^m_s(\vec{\Phi}(z))}{s^2}\ dl_{\p B_r}
\end{array}
\ee
Using (\ref{reg-III.5-a}) we then obtain
\be
\label{reg-III.14}
\begin{array}{l}
\ds\limsup_{s\rightarrow 0}\lf|\int_{B_r(x)\cap{\mathcal B}_\ep} \xi_\ep(\vec{\Phi})\, \nabla(\eta_s(\vec{\Phi}))\,F(\vec{\Phi})_{ij}\, \nabla\vec{\Phi}_i\ dy^2\rg| \\[5mm]
\ds\quad\le  C\ \|F\|_\infty\ \int_{\p B_r(x)\cap {\mathcal G}_\ep}\ |\nabla\vec{\Phi}|(z)\ \theta^2_\ep(\vec{\Phi}(y))\ dl_{\p B_r}(z)
\end{array}
\ee
and using the uniform bound on the density which itself comes from the monotonicity formula
\be
\label{reg-III.15}
\begin{array}{l}
\ds\limsup_{s\rightarrow 0}\lf|\int_{B_r(x)\cap{\mathcal B}_\ep} \xi_\ep(\vec{\Phi})\, \nabla(\eta_s(\vec{\Phi}))\,F(\vec{\Phi})_{ij}\, \nabla\vec{\Phi}_i\ dy^2\rg| \\[5mm]
\ds\ \le \ C\ \|F\|_\infty\ \int_{\p B_r(x)\cap {\mathcal G}_\ep\cap \vec{\Phi}^{-1}({\mathfrak B}_\ep)}\ |\nabla\vec{\Phi}|(z)\ dl_{\p B_r}(z)\le \ C\ \|F\|_\infty\ \delta
\end{array}
\ee
where we have used (\ref{reg-III.7-a}). 
Observe that since we have cut-off $\Gamma_{r,x}\setminus\Gamma_{r,x}^\ep$ by multipying by $\xi_\ep(\vec{\Phi})$ we have for $s$ small enough
\be
\label{reg-III.16}
\begin{array}{l}
\ds\lim_{s\rightarrow 0}\int_{B_r(x)\cap {\mathcal G}_\ep} \xi_\ep(\vec{\Phi})\, \nabla(\eta_s(\vec{\Phi}))\,F(\vec{\Phi})_{ij}\, \nabla\vec{\Phi}_i\ dy^2\\[5mm]
\ds\lim_{s\rightarrow 0}\int_{B_r(x)\cap {\mathcal G}_\ep}\xi_\ep(\vec{\Phi})\, \nabla\lf(\eta\lf(\frac{\mbox{dist}(\vec{\Phi}(y),\Gamma_{r,x}^\ep)}{s}\rg)\rg)\,F(\vec{\Phi})_{ij}\, \nabla\vec{\Phi}_i\ dy^2
\end{array}
\ee
Using the coarea formula this gives
\be
\label{reg-III.17}
\begin{array}{l}
\ds\lim_{s\rightarrow 0}\sum_{i=1}^m\int_{B_r(x)\cap {\mathcal G}_\ep}\xi_\ep(\vec{\Phi})\, \nabla\lf(\eta\lf(\frac{\mbox{dist}(\vec{\Phi}(y),\Gamma_{r,x}^\ep)}{s}\rg)\rg)\,F(\vec{\Phi})_{ij}\, \nabla\vec{\Phi}_i\ dy^2\\[5mm]
\ds\quad=\lim_{s\rightarrow 0}\frac{1}{s}\int_0^{2s}\chi'\lf(\frac{\sigma}{s}\rg)\ d\sigma\int_{\mbox{dist}(\vec{p},\Gamma_{r,x}^\ep)=\sigma}\xi_\ep(\vec{p})\ F_j(\vec{p})\cdot\nu\ d{\mathcal H}^1\res\Sigma_{r,x}^\ep
\end{array}
\ee
where $\Sigma_{r,x}^\ep$ is the immersed sub-manifold $\vec{\Phi}(B_r(x))\cap{\mathcal G}_\ep$ and $\vec{\nu}$ is the unit exterior vector in this
sub-manifold orthogonal to the level set $\mbox{dist}(\vec{p},\Gamma_{r,x}^\ep)=\sigma$ and $F_j(\vec{p})\cdot\nu=\sum_{i=1}^mF_{ji}(\vec{p})\,\nu_i$. Since $\vec{\Phi}$ is an immersion in a neighborhood of $\p B_r\cap {\mathcal G}_\ep$, for $\sigma$ small enough and being a regular value of
$f_{r,x}^\ep(\vec{p}):=\mbox{dist}(\vec{p},\Gamma_{r,x}^\ep)$ the level set is made of the following union
\[
\lf(f_{r,x}^\ep\rg)^{-1}(\sigma)=\gamma_{r,x}^\ep(\sigma)\cup_{\al\in A_\sigma} \p \om_\al(\sigma)
\]
where $\gamma_{r,x}^\ep(\sigma)$ is a smooth curve converging to $\Gamma_{r,x}^\ep$ and   $\om_\al(\sigma)$ are subdomains of $B_r(x)$ included in $\mbox{dist}(\vec{\Phi}(y),\Gamma_{r,x}^\ep)^{-1}([0,\sigma])$. The Taylor expansion
of $\vec{\Phi}$ with respect to each point $y\in \p B_r(x)$ gives
\be
\label{reg-III.18}
\begin{array}{l}
\ds\lim_{\sigma\rightarrow 0}\int_{\gamma_{r,x}^\ep(\sigma)}\xi_\ep(\vec{p})\ F_j(\vec{p})\cdot\nu\ d{\mathcal H}^1\res\Sigma_{r,x}^\ep=
\int_{\p B_r(x)\cap{\mathcal G}_\ep}\xi_\ep(\vec{\Phi})\ \sum_{i=1}^mF_{ij}(\vec{\Phi})\, \frac{\p\vec{\Phi}_i}{\p r}\ dl_{\p B_r}
\end{array}
\ee
Hence
\be
\label{reg-III.19}
\begin{array}{l}
\ds\lim_{s\rightarrow 0}\frac{1}{s}\int_0^{2s}\chi'\lf(\frac{\sigma}{s}\rg)\ d\sigma\int_{\gamma_{r,x}^\ep(\sigma)}\xi_\ep(\vec{p})\ F_j(\vec{p})\cdot\nu\ d{\mathcal H}^1\res\Sigma_{r,x}^\ep\\[5mm]
\ds\quad=\int_{\p B_r(x)\cap{\mathcal G}_\ep}\xi_\ep(\vec{\Phi})\ \sum_{i=1}^mF_{ij}(\vec{\Phi})\, \frac{\p\vec{\Phi}_i}{\p r}\ dl_{\p B_r}
\end{array}
\ee
For the other contributions we have
\be
\label{reg-III.20}
\begin{array}{l}
\ds\lim_{s\rightarrow 0}\frac{1}{s}\int_0^{2s}\chi'\lf(\frac{\sigma}{s}\rg)\ d\sigma\sum_{\al\in A_\sigma}\int_{\p\om_\al(\sigma)}\xi_\ep(\vec{p})\ F_j(\vec{p})\cdot\nu\ d{\mathcal H}^1\res\Sigma_{r,x}^\ep\\[5mm]
\ds\quad\lim_{s\rightarrow 0}\frac{1}{s}\int_0^{2s}\chi'\lf(\frac{\sigma}{s}\rg)\ d\sigma\sum_{\al\in A_\sigma}\int_{\om_\al(\sigma)}\mbox{div}_{\Sigma_{r,x}^\ep}\lf(\xi_\ep(\vec{p})\ F_j(\vec{p})\rg)\ d{\mathcal H}^2\res\Sigma_{r,x}^\ep
\end{array}
\ee
Since $\om_\al(\sigma)$ is included in a $\sigma$ neighborhood of the smooth curve $\Gamma_{r,x}^\ep$, using the monotonicity formula, for $\sigma$ small enough, covering such a neighborhood by $\simeq\sigma^{-1}\ {\mathcal H}^1(\Gamma_{r,x}^\ep)$ balls of radius $2\sigma$ we have the following bound
\be
\label{reg-III.21}
\sum_{\al\in A_\sigma}{\mathcal H}^2(\om_\al(\sigma))\le C\ \sigma\ {\mathcal H}^1(\Gamma_{r,x}^\ep)
 \sup_{t,\vec{p}}\frac{{\mathcal H}^2(\Sigma_{r,x}^\ep\cap B^m_t(\vec{p}))}{t^2} \le C_\ep\ \sigma
\ee
Hence we have
\be
\label{reg-III.22}
\begin{array}{l}
\ds\lim_{s\rightarrow 0}\frac{1}{s}\int_0^{2s}\chi'\lf(\frac{\sigma}{s}\rg)\ d\sigma\sum_{\al\in A_\sigma}\int_{\p\om_\al(\sigma)}\xi_\ep(\vec{p})\ F_j(\vec{p})\cdot\nu\ d{\mathcal H}^1\res\Sigma_{r,x}^\ep\\[5mm]
\ds\quad\le \ds\lim_{s\rightarrow 0}\frac{C_\ep}{s}\int_0^{2s}\|\chi'\|_\infty\ \|\mbox{div}_{\Sigma_{r,x}^\ep}\lf(\xi_\ep(\vec{p})\ F_j(\vec{p})\rg)\|_{L^\infty(\Sigma_{r,x}^\ep)}\ \sigma\  \ d\sigma\le C_\ep\ \sigma
\end{array}
\ee
Combining (\ref{reg-III.16})... (\ref{reg-III.22}) we obtain
\be
\label{reg-III.23}
\begin{array}{l}
\ds\lim_{s\rightarrow 0}\int_{B_r(x)\cap {\mathcal G}_\ep} \xi_\ep(\vec{\Phi})\, \nabla(\eta_s(\vec{\Phi}))\,F(\vec{\Phi})_{ij}\, \nabla\vec{\Phi}_i\ dy^2\\[5mm]
\ds\quad=\int_{\p B_r(x)\cap{\mathcal G}_\ep}\xi_\ep(\vec{\Phi})\ \sum_{i=1}^mF_{ij}(\vec{\Phi})\, \frac{\p\vec{\Phi}_i}{\p r}\ dl_{\p B_r}
\end{array}
\ee
Collecting (\ref{reg-III.7-a}), (\ref{reg-III.9}), (\ref{reg-III.15}) and (\ref{reg-III.23}) we obtain
\be
\label{reg-III.24}
\begin{array}{l}
\ds \lf|\sum_{i=1}^m\int_{\p B_r(x)} F(\vec{\Phi}){ij}\,\frac{\p \vec{\Phi}_i}{\p r}\, dl_{\p B_r}-\int_{B_r(x)} \nabla(F(\vec{\Phi}){ij})\, \nabla\vec{\Phi}_i\ dy^2 \rg.\\[5mm]
\ds\quad\lf.+\int_{B_r(x)}F(\vec{\Phi})_{ij}\, A_i(\vec{\Phi})(d\vec{\Phi},d\vec{\Phi})\ dy^2\rg|<C(\|F\|_\infty)\ \delta
\end{array}
\ee
This holds for any $\delta>0$ and hence lemma~\ref{lm-reg-III.5} is proved. \hfill $\Box$

\subsection{Proof of theorem~\ref{th-I.1}.} Let $\varphi\in C^\infty_0(D^2)$, be a non negative function. The co-area formula gives
\be
\label{reg-III.25}
\int_{D^2}\Delta\varphi\ \vec{\Phi}\ dx^2=-\int_{D^2}\nabla\varphi\cdot\nabla\vec{\Phi}=\int_{0}^{+\infty}dt\int_{\varphi^{-1}\{t\}}\p_\nu\vec{\Phi}\ dl_{\varphi^{-1}\{t\}}
\ee
Using lemma~\ref{lm-reg-III.5} for $F_{ij}:=\delta_{ij}$, we obtain
\be
\label{reg-III.26}
\begin{array}{l}
\ds\int_{D^2}\Delta\varphi\ \vec{\Phi}\ dx^2=-\int_0^{+\infty}dt\, \int_{\varphi^{-1}((t,+\infty))} A(\vec{\Phi})(d\vec{\Phi},d\vec{\Phi})\ dy^2\\[5mm]
\ds\quad=-\int_0^{+\infty}\int_{D^2}{\mathbf 1}_{\varphi(x)>t} \ A(\vec{\Phi})(d\vec{\Phi},d\vec{\Phi})\ dy^2\ dt=-\int_{D^2}\varphi(x)\ A(\vec{\Phi})(d\vec{\Phi},d\vec{\Phi})\ dy^2
\end{array}
\ee
This implies that $\vec{\Phi}$ satisfies weakly the harmonic map equation into $N^n$ 
\[
-\Delta\vec{\Phi}=A(\vec{\Phi})(d\vec{\Phi},d\vec{\Phi})\quad\quad\mbox{ in }{\mathcal D}'(D^2)
\]
and using H\'elein's regularity result \cite{Hel}, we prove theorem~\ref{th-I.1}.\hfill $\Box$
 \renewcommand{\theequation}{A.\arabic{equation}}
\renewcommand{\theTh}{A.\arabic{Th}}
\renewcommand{\theProp}{A.\arabic{Prop}}
\renewcommand{\theLma}{A.\arabic{Lma}}
\renewcommand{\theCo}{A.\arabic{Co}}
\renewcommand{\theRm}{A.\arabic{Rm}}
\renewcommand{\theequation}{A.\arabic{equation}}
\setcounter{equation}{0} 
\reset
\appendix
\section{Appendix}

\begin{Prop}
\label{prop-A-1} Let $N^n$ be a $C^2$ sub-manifold of the euclidian space ${\R}^m$. Let $(\Sigma,h)$ be a compact Riemann surface (equipped with a metric compatible with the complex structure)
possibly with boundary.
Let $\vec{\Phi}$ be a map in $W^{1,2}(\Sigma,N^n)$. Assume $\vec{\Phi}$ is weakly conformal and continuous on $\p\Sigma$. The integer rectifiable varifold associated to $(\vec{\Phi},\Sigma)$ is stationary in $N^n\setminus\vec{\Phi}(\p\Sigma)$ if and only if
\[
\forall\ {F}\in C^\infty_0(N^n\setminus\vec{\Phi}(\p\Sigma),{\R}^m)\quad\int_{\Sigma} \lf[\lf<d(F(\vec{\Phi})), d\vec{\Phi}\rg>_h- F(\vec{\Phi})\ A(\vec{\Phi})(d\vec{\Phi},d\vec{\Phi})_h\rg]\ dvol_h=0
\]
where $A(\vec{q})(\vec{X},\vec{Y})$ denotes the second fundamental form of $N^n$ at the point $\vec{q}$ and acting on the pair of vectors $(\vec{X},\vec{Y})$ and by an abuse of notation we write
\[
\ A(\vec{\Phi})(d\vec{\Phi},d\vec{\Phi})_h:=\sum_{i,j=1}^2h_{ij}\ A(\vec{\Phi})(\p_{x_i}\vec{\Phi},\p_{x_j}\vec{\Phi})\quad.
\].\hfill $\Box$
\end{Prop}
\noindent{\bf Proof of proposition~\ref{prop-A-1}.}
For any $\vec{q}\in N^n$ one denotes by $P_T(\vec{q})$ the symmetric matrix giving the orthogonal projection onto $T_{\vec{q}}N^n$. 
The integer rectifiable varifold given by $(\vec{\Phi},\Sigma)$  is by defintion the following {\it Radon measure} on $G_2(T{\R}^m)$ the Grassman bundle of un-oriented 2-planes over ${\R}^m$ given by
\[
\forall \ \phi\in C^\infty(G_2(T{\R}^m))\quad\quad{\mathbf v}_{\vec{\Phi}}(\phi)=\int_{G_2(T{\R}^m)}\phi(S,\vec{q})\ dV_{\vec{\Phi}}(S,\vec{q}):=\int_{\Sigma} \phi(\vec{\Phi}_\ast(T_{x}\Sigma),{\vec{\Phi}(x)})\ dvol_{g_{\vec{\Phi}}}
\]  
By definition (see \cite{All}), the varifold ${\mathbf v}_{\vec{\Phi}}$ is stationary in $N^n$ if
\be
\label{A-1}
\forall\ {F}\in C^\infty_0(N^n\setminus\vec{\Phi}(\p\Sigma),{\R}^m)\quad\int_{\Sigma} \mbox{div}_S(P_T\, F)(\vec{q})\ dV_{\vec{\Phi}}(S,\vec{q})=0
\ee
In local conformal coordinates at a point where $|\p_{x_1}\vec{\Phi}|=|\p_{x_2}\vec{\Phi}|=e^\la$, introducing the orthonormal basis of $S:=\vec{\Phi}_\ast T_x\Sigma$ given by $\vec{e}_i:=e^{-\la}\p_{x_i}\vec{\Phi}$, one has by definition
\[
\mbox{div}_{\vec{\Phi}_\ast T_x\Sigma}(P_T\, F)(\vec{\Phi}):=\sum_{i=1}^2\p_{\vec{e}_i}(P_T\, F)(\vec{\Phi})\cdot\vec{e}_i=\sum_{i=1}^2\sum_{k=1}^me_i^k\ \p_{z_k}(P_T\, F)(\vec{\Phi})\cdot\vec{e}_i
\]
where $\vec{e}_i:=\sum_{k=1}^me_i^k\ \p_{z_k}$. Hence we have
\[
\begin{array}{l}
\ds\mbox{div}_{\vec{\Phi}_\ast T_x\Sigma}(P_T\, F)(\vec{\Phi})=e^{-2\la}\,\sum_{i=1}^2\sum_{k=1}^m\p_{x_i}\Phi^k\ \p_{z_k}(P_T\, F)(\vec{\Phi})\cdot\p_{x_i}\vec{\Phi}\ \\[5mm]
\ds\quad=e^{-2\la}\,\nabla (F(\vec{\Phi}))\cdot\nabla\vec{\Phi}- F(\vec{\Phi})\ \sum_{i=1}^2 e^{-2\la}\ A(\vec{\Phi})(\p_{x_i}\vec{\Phi},\p_{x_i}\vec{\Phi})
\end{array}
\]
where we used respectively that $P_T(\vec{\Phi})\nabla\vec{\Phi}=\nabla\vec{\Phi}$ and that $A(\vec{q})(\vec{X},\vec{Y})=-\,\p_{\vec{X}}P_T(\vec{q})\cdot\vec{Y}$.  Multiplying by $dvol_{g_{\vec{\Phi}}}= e^{2\la}\ dx_1\wedge dx_2$ we obtain at almost every point $x$ where $\nabla\vec{\Phi}(x)\ne0$
\be
\label{A-2}
\mbox{div}_{\vec{\Phi}_\ast T_x\Sigma}(P_T\, F)(\vec{\Phi})\ dvol_{g_{\vec{\Phi}}}= \lf[\lf<d(F(\vec{\Phi})), d\vec{\Phi}\rg>_{g_{\vec{\Phi}}}- F(\vec{\Phi})\ A(\vec{\Phi})(d\vec{\Phi},d\vec{\Phi})_{g_{\vec{\Phi}}}\rg]\ dvol_{g_{\vec{\Phi}}}
\ee
 This concludes the proof
of the lemma.\hfill $\Box$

\end{document}